\documentclass[submission]{FPSAC2021}

\usepackage{mathrsfs}
\usepackage{stackengine}
\usepackage{todonotes}

\newtheorem{thm}{Theorem}

\theoremstyle{definition}
\newtheorem{defn}{Definition}
\newtheorem{algorithm}{Algorithm}

\newtheorem{example}{Example}

\usepackage{lipsum}
\usepackage{ulem}

\usepackage{wrapfig}

\usepackage{color, tikz}
\usetikzlibrary{positioning}
\usetikzlibrary{arrows}

\usetikzlibrary{shapes}

\title[]{The combinatorial PT-DT correspondence}

\author[Helen Jenne, Gautam Webb, and Benjamin Young]{Helen Jenne\thanks{\href{mailto:helen.jenne@lmpt.univ-tours.fr}{helen.jenne@lmpt.univ-tours.fr}. Helen Jenne has received funding from the European Research Council (ERC) through the European Union's Horizon 2020 research and innovation programme under the Grant Agreement No 759702. }\addressmark{1}, Gautam Webb\thanks{\href{mailto:gwebb@uoregon.edu}{gwebb@uoregon.edu} }\addressmark{2}, \and Benjamin Young \thanks{\href{mailto:bjy@uoregon.edu}{bjy@uoregon.edu}. Benjamin Young was partially supported by the Knut and Alice Wallenberg Foundation Grant KAW:2010.0063.}\addressmark{2} }

\address{\addressmark{1}CNRS, Institut Denis Poisson, Universit\'e de Tours and Universit\'e d'Orl\'eans, France \\ \addressmark{2}Department of Mathematics, University of Oregon, Eugene, OR, USA}

\received{\today}

\abstract{We resolve an open conjecture from algebraic geometry, which states that two generating functions for plane partition-like objects (the "box-counting" formulae for the Calabi-Yau topological vertices in Donaldson-Thomas theory and Pandharipande-Thomas theory) are equal up to a factor of MacMahon's generating function for plane partitions. The main tools in our proof are a Desnanot-Jacobi-type "condensation" identity, and a novel application of the tripartite double-dimer model of Kenyon-Wilson.}


\begin{document}

\newcommand{\sA}{\mathscr{A}}
\newcommand{\sB}{\mathscr{B}}
\newcommand{\sAB}{\mathscr{A}\!\!\mathscr{B}}
\newcommand{\cAB}{\mathcal{A}\mathcal{B}}
\newcommand{\Cyl}{\text{Cyl}}
\newcommand{\I}{\text{I}}
\newcommand{\II}{\text{I}\!\protect\rule{0.015in}{0in}\text{I}}
\newcommand{\III}{\text{I}\!\protect\rule{0.015in}{0in}\text{I}\!\protect\rule{0.015in}{0in}\text{I}}
\newcommand{\base}{\text{base}}

\maketitle

\section{Introduction}

Donaldson-Thomas (DT) theory and Pandharipande-Thomas (PT) theory are branches of enumerative geometry closely related to mirror symmetry and string theory. In both theories, generating functions arise  known as the \emph{combinatorial Calabi-Yau topological vertices}. 
These generating functions enumerate seemingly different plane partition-like objects. In this paper, we prove that the generating functions coincide up to a factor of $M(q)$, MacMahon's generating function for plane partitions \cite{macmahon}. Our result, taken together with a substantial body of geometric work, proves a geometric conjecture in the foundational work of Pandharipande-Thomas theory which has been open for over 10 years.

The generating function from Donaldson-Thomas theory is known as the DT topological vertex. Denoted $V(\mu_1, \mu_2, \mu_3)$, where each $\mu_i$ is a partition, it counts plane partitions {\em asymptotic to $\mu_1, \mu_2$, and $\mu_3$} (see Section~\ref{sec:DTboxconfigs}). 
The PT topological vertex, denoted by $W(\mu_1, \mu_2, \mu_3)$, is a generating function for a certain class of finitely generated $\mathbb{C}[x_1,x_2,x_3]$-modules (see Section~\ref{sec:PTboxconfigs}).

We prove that

\begin{thm}\cite[Calabi-Yau case of Conjecture 4]{PT2}
\label{thm:ptdt}
\begin{equation}
	\label{eqn:ptdt}
	V(\mu_1, \mu_2, \mu_3) = M(q) W(\mu_1, \mu_2, \mu_3),
\end{equation}
where $M(q) = \prod\limits_{i\geq 1}\left(1-q^i\right)^{-i}$.
\end{thm}

The geometric corollary of this theorem is a proof of
Theorem/Conjecture 2 of \cite{PT2}, which, loosely speaking, states that $W(\mu_1, \mu_2, \mu_3)$ computes the local contribution to the Pandharipande-Thomas generating function. The proof of this corollary combines Theorem~\ref{thm:ptdt} with the analogous result in DT theory \cite{mnop1, mnop2, mpt} along with \cite[Section 4.1.2]{moop}; it is a consequence of the fact that both DT and PT theory give the same invariants as a third theory, Gromov-Witten theory.\footnote{In~\cite{PT2, mnop1, mnop2} and in general elsewhere in the geometry literature, all of the formulas have $q$ replaced by $-q$.  The sign is there for geometric reasons which are immaterial for us.}

The combinatorics problems which we solve are stated in the geometry literature as ``box-counting'' problems; that is, the objects of interest are plane partition-like.  The following bijections are well-known:
\[
	\text{\stackanchor{dimer configurations on}{the honeycomb graph}}
	\leftrightarrow
	\text{plane partitions}
	\leftrightarrow
	\text{\stackanchor{finite-length monomial}{ideals in $\mathbb{C}[x_1,x_2,x_3]$}}
\]
The first one is a 3D version of the correspondence between a partition and its Maya diagram; it is stated explicitly in Section~\ref{sec:DTtheoryAndDimers}. 
We use essentially the same correspondence to give a dimer description of the DT topological vertex $V(\mu_1, \mu_2, \mu_3)$.
On the PT side, the correspondences are
\[
	\text{\stackanchor{tripartite double-dimer configs.}{on the honeycomb graph}}
	\stackrel{(1)}{\leftrightarrow}
	\text{\stackanchor{labelled box}{configurations}}
	\stackrel{(2)}{\leftrightarrow}
	\text{\stackanchor{$\mathbb{C}[x_1, x_2, x_3]$-modules}{$(M_1 \oplus M_2 \oplus M_3)/\left<(1,1,1)\right>$}}
	\]
The correspondence (1) is new, as far as we are aware.  We describe labelled box configurations, and the generating functions for them which arise in PT theory, carefully in Section~\ref{sec:PT}.  Interestingly, though it is a purely combinatorial correspondence, it is not bijective---rather, it is a weight-preserving, 1-to-many correspondence.
Here $M_1 \subseteq \mathbb{C}[x_1, x_1^{-1}, x_2, x_3]$ is spanned by all monomials $x_1^ix_2^jx_3^k$ where $i \in \mathbb{Z}$ and $(j,k)$ ranges over some fixed partition $\mu_1$, with $M_2, M_3$ defined similarly; the quotient is killing the diagonal of the direct sum.  

The correspondence (2) is incidental to this work and is described in~\cite{PT2}; nor will we need to discuss the structure of the modules in the codomain. We expect that our methods will be relevant in other similar situations (one such situation arises in Rank 2 DT theory~\cite{GKY2017}) and we would be eager to learn of other instances in which our techniques would apply.

We prove Theorem~\ref{thm:ptdt} by observing that both $V(\mu_1, \mu_2, \mu_3)/M(q)$ and $W(\mu_1, \mu_2, \mu_3)$ are the unique solution of the same recurrence, with the same initial conditions.  The recurrence in question is called the \emph{condensation recurrence}; we postpone its definition to Section~\ref{sec:definitions}, after we have made the required definitions.

Viewed as a recurrence in $\mu_1$ and $\mu_2$, Equation~\eqref{eqn:vertex_condensation} uniquely characterizes $V$ and $W$. The base case is when one of the three partitions $\mu_i$ is equal to zero; Equation \eqref{eqn:ptdt} is known to hold in this situation~\cite{PT2}.

When recast in terms of the dimer model, $V(\mu_1, \mu_2, \mu_3)$ is easily seen to satisfy Equation~\eqref{eqn:vertex_condensation} by Kuo's \emph{graphical condensation}~\cite{kuo}; this is essentially the content of Section~\ref{sec:DT}. 

Showing that $W(\mu_1, \mu_2, \mu_3)$ satisfies Equation~\eqref{eqn:vertex_condensation} is considerably more intricate, but once we translate to the double-dimer model, the bulk of the work was done elsewhere, in work of Jenne~\cite{jenne}.  Essentially,~\cite{jenne} evaluates a certain determinant by the classical Desnanot-Jacobi identity, and then interprets all six terms in the identity in terms of $W$. 

The full version of this abstract will appear in \cite{JWY}; proofs have been omitted due to space constraints.

\section{Definitions}
\label{sec:definitions}
Fix partitions $\mu_1, \mu_2, \mu_3$.  For this paper, we identify $\mu_i$ with the coordinates of the boxes of its Young diagram, with the corner of the diagram located at $(0,0)$.  Define the following subsets of $\mathbb{Z}^3$, thought of as sets of boxes: $\Cyl_1 = \{(x,u,v) \in \mathbb{Z}^3 \;|\; (u,v) \in \mu_1\}, \Cyl_2 = \{(v,y,u) \in \mathbb{Z}^3 \;|\; (u,v) \in \mu_2\}$, and $\Cyl_3 = \{(u,v,z) \in \mathbb{Z}^3 \;|\; (u,v) \in \mu_3\}.$

Moreover, let $\mathbb{Z}^3_{\geq 0}$ denote the integer points in the first octant (including the coordinate planes and axes). Let $\Cyl_i^{+} = \Cyl_i \cap \mathbb{Z}^3_{\geq 0}$ and $\Cyl_i^{-} = \Cyl_i \setminus \mathbb{Z}^3_{\geq 0}$.
Finally, let
\begin{align*}
	&&\II_{\bar{1}} &= \Cyl_2 \cap \Cyl_3\setminus\Cyl_1, \\
	\I^- &= \Cyl_1^- \cup \Cyl_2^- \cup \Cyl_3^-, &
	\II_{\bar{2}} &= \Cyl_3 \cap \Cyl_1\setminus\Cyl_2, &
	\II &= \II_{\bar{1}} \cup \II_{\bar{2}} \cup \II_{\bar{3}}, \\
	\I^+ &= \Cyl_1^+ \cup \Cyl_2^+ \cup \Cyl_3^+, &
	\II_{\bar{3}} &= \Cyl_1 \cap \Cyl_2\setminus\Cyl_3, & \III &= \Cyl_1 \cap \Cyl_2 \cap \Cyl_3.
\end{align*}

We will need the following standard notions of Maya diagrams.
If $\lambda = (\lambda_1, \lambda_2, \ldots, \lambda_k)$ is a partition with $k$ parts, define $\lambda_t=0$ for $t>k$. The \emph{Maya diagram of $\lambda$} is the set $\{\lambda_t - t + \frac{1}{2}\} \subseteq \mathbb{Z}+\frac{1}{2}$.  We frequently associate a partition with its Maya diagram by drawing a Maya diagram as a doubly infinite sequence of beads and holes, indexed by $\mathbb{Z}+\frac{1}{2}$, with the beads representing elements of the above set.  For instance, the Maya diagrams of the empty partition and of the partition $\lambda = (4,2,1)$ are the sets $\{ -\frac{1}{2}, -\frac{3}{2}, \ldots \}$ and $\{ \frac{7}{2}, \frac{1}{2}, -\frac{3}{2},-\frac{7}{2}, -\frac{9}{2},  \ldots \}$, respectively, which are drawn as
\[
\cdots  \circ \circ \circ | \bullet \bullet \bullet \cdots
\qquad
\text{and}
\qquad
\cdots \circ \circ \circ \bullet \circ \circ \bullet | \circ \bullet \circ \bullet \bullet \bullet \cdots. 
\]
When convenient, we simply mark the location of 0 with a vertical line, rather than labelling the beads with elements of $\mathbb{Z}+\frac{1}{2}$.
Conversely, if $S$ is a subset of $\mathbb{Z}+\frac{1}{2}$,  define $S^+ = \{x \in S \;|\; x > 0\}$ and $S^- = \{x \in \mathbb{Z}+\frac{1}{2} \setminus S \;|\; x < 0\}$.  If both $S^+$ and $S^-$ are finite, then define the \emph{charge} of $S$, $c(S)$, to be $|S^+|-|S^-|$; then it is easy to check that the set $\{s - c(S) \;:\;s \in S\}$ is the Maya diagram of some partition $\lambda$; we say that $S$ itself is the \emph{charge $c(S)$ Maya diagram of $\lambda$}.

If $\lambda$ is a partition with Maya diagram $S$, let
$\lambda^r$ (resp. $\lambda^c$) be the partition associated to the charge $-1$ (resp. $1$) Maya diagram $S \setminus \{\min S^+\}$ (resp. $S \cup \{ \max S^-\}$). Let $\lambda^{rc}$ be the partition associated to the Maya diagram $(S\setminus\{ \min{S^+}\})\cup\{\max{S^-}\}$.


\begin{wrapfigure}{r}{.6\textwidth}
\centering
\includegraphics[width=4in]{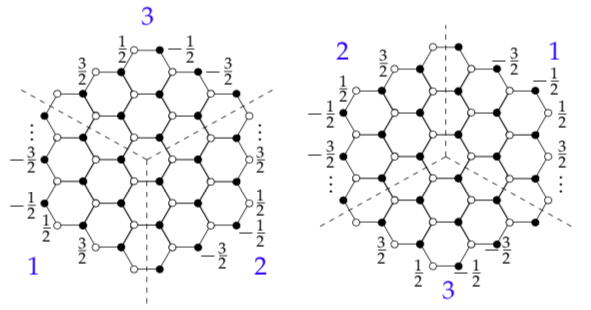}
\caption{The division of $H(3)$ into sectors for DT and PT}
\label{fig:sectors}
\end{wrapfigure}
In both PT and DT, it will be convenient to divide the $N \times N \times N$ honeycomb graph $H(N)$ into three sectors and label some of the vertices on the outer face as shown in Figure~\ref{fig:sectors} for $H(3)$. We remark that the division into sectors makes sense as $N \to \infty$.  The reason for these particular choices of labels is that we will need to specify these specific vertices, both in DT and PT, based on the Maya diagrams of various partitions.

Finally, let $X=X(\mu_1, \mu_2, \mu_3)$ be a power series in $q$, depending on three partitions $\mu_1, \mu_2, \mu_3$, which is symmetric with respect to cyclic permutation of these partitions.  We shall be interested in solutions $X$ to the following functional equation:
\begin{equation}\small
	\label{eqn:vertex_condensation}
	q^A
	X(\mu_1, \mu_2, \mu_3)
	X(\mu_1^{rc}, \mu_2^{rc}, \mu_3)
	=
	q^B
	X(\mu_1^{rc}, \mu_2, \mu_3)
	X(\mu_1, \mu_2^{rc}, \mu_3)
	+
	q^C
	X(\mu_1^{r}, \mu_2^{c}, \mu_3)
	X(\mu_1^{c}, \mu_2^{r}, \mu_3).
\end{equation}
Here, $A,B,C$ are certain explicit constants depending on $\mu_1, \mu_2, \mu_3$ which we don't define in this extended abstract. These constants are discussed further in Section~\ref{sec:DTcond}.~


Since the partitions $\mu_i^r$, $\mu_i^c$, $\mu_i^{rc}$ are all smaller, in some sense, than $\mu_i$,
and since none of the topological vertex terms are equal to zero, we can divide both sides of the condensation recurrence by $V(\mu_1^{rc}, \mu_2^{rc}, \mu_3)$ and obtain a recursive characterization of $V(\mu_1, \mu_2, \mu_3)$. Note also that $V(\mu_1, \mu_2, \mu_3) = V(\mu_2, \mu_3, \mu_1)$ by symmetry - so we can say that $V(\mu_1, \mu_2, \mu_3)$ is the \emph{unique} power series which satisfies the condensation recurrence, where we take the base cases to be the (known) value of $V(\mu_1, \mu_2, \emptyset)$ for all partitions $\mu_1, \mu_2$.

\section{DT}
\label{sec:DT}
\subsection{DT box configurations}
\label{sec:DTboxconfigs}

We say that a \emph{plane partition asymptotic to ($\mu_1, \mu_2, \mu_3$)} is an order ideal under the product order in $\mathbb{Z}^3_{\geq 0}$ which contains $\I^+ \cup \II \cup \III$, together with only finitely many other points in $\mathbb{Z}_{\geq 0}^3$.  We let $P(\mu_1,\mu_2,\mu_3)$ denote the set of plane partitions asymptotic to $\mu_1, \mu_2, \mu_3$.  

If any of $\mu_1, \mu_2, \mu_3$ are nonzero, then every $\pi \in P(\mu_1, \mu_2, \mu_3)$ is an infinite subset of $\mathbb{Z}^3_{\geq 0}$.  We define
$w(\pi) = |\pi \setminus (\I^+ \cup \II \cup \III)| - |\II| -2|\III|, $
the customary measure of ``size'' of such a plane partition in the geometry literature (see, for instance,~\cite{mnop1}).   

Define 
	\[ V(\mu_1, \mu_2, \mu_3) = \sum\limits_{\pi \in P(\mu_1, \mu_2, \mu_3)} q^{w(\pi)}.\]
We call $V(\mu_1, \mu_2, \mu_3)$ the \emph{topological vertex in Donaldson-Thomas theory}.  Note that if $\pi \in P(\emptyset, \emptyset, \emptyset)$ with $|\pi|=n$, then $\pi$ is a plane partition of $n$ in the conventional sense, that is, a finite array of integers such that each row and column is a weakly decreasing sequence of nonnegative integers. Thus MacMahon's enumeration of plane partitions~\cite{macmahon} gives us
	$V(\emptyset, \emptyset, \emptyset) = \prod_{i = 1}^{\infty}\left(1-q^i\right)^{-i}$.

In~\cite{orv}, there is an expansion of $V(\mu_1,\mu_2,\mu_3)$ in terms of Schur functions.  However, since no similar expansion is known in PT theory, this expansion 
does not help prove Theorem~\ref{thm:ptdt}. 
\subsection{DT theory and the dimer model}
\label{sec:DTtheoryAndDimers}

Before giving the dimer description of $V(\mu_1,\mu_2,\mu_3)$, we review the correspondence between plane partitions and dimer configurations of a honeycomb graph. By representing each integer $i$ in a plane partition as a stack of $i$ unit boxes, a plane partition can be visualized as a collection of boxes which is stacked stably in the positive octant, with gravity pulling them in the direction $(-1, -1, -1)$.
This collection of boxes can be viewed as a lozenge tiling of a hexagonal region of triangles, which is equivalent to a dimer configuration (also called a perfect matching) of its dual graph.  

\begin{figure}[htb]
\begin{center}
\includegraphics[width=6in]{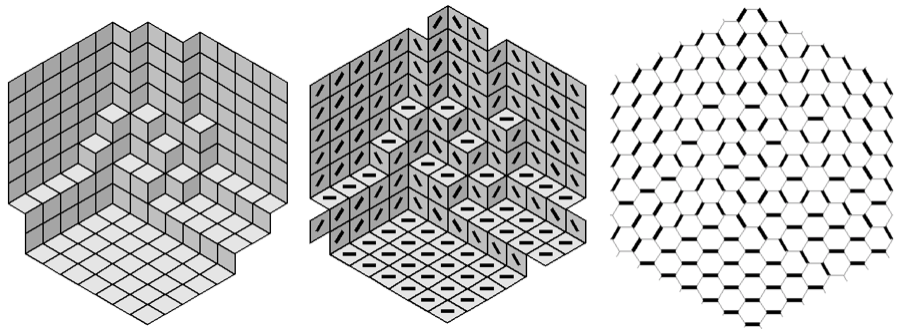}
\end{center}
\caption{Shown left is a plane partition $\pi$ asymptotic to $(\mu_1, \mu_2, \mu_3)$, where $\mu_1 =(1, 1)$, $\mu_2 = \mu_3 = (2, 1, 1)$, $|\II| = 9$, $|\III| = 3$, and $w(\pi) =13-|\II|-2|\III| = -2$. We see that $\pi$ is equivalent to a tiling, which is truncated in the center image so that it corresponds to 
 a dimer configuration of $H(7)$ with a few vertices on the outer face deleted. }
\label{fig:DTdimers}
\end{figure}

Just as a plane partition can be visualized as a collection of boxes, a plane partition asymptotic to $(\mu_1, \mu_2, \mu_3)$ can be visualized as a collection of boxes, as shown in Figure~\ref{fig:DTdimers}, left picture.  Moreover, a version of the above correspondence puts these box collections in bijection with dimer configurations on the honeycomb graph $H(N)$ with some outer vertices removed, which we call $H(N; \mu_1, \mu_2, \mu_3)$. Specifically, $H(N; \mu_1, \mu_2, \mu_3)$ is constructed as follows. Let $S_i$ be the Maya diagram of $\mu_i$. Construct the sets $S_i^+$, $S_i^-$ for $i=1,2,3$ and then remove the vertices with the labels in $S_i^+\cup S_i^-$ from sector $i$ of $H(N)$ (here, we are referring to the labelling of the boundary vertices illustrated in Figure~\ref{fig:sectors}, left picture).

\subsection{The condensation recurrence in DT theory}
\label{sec:DTcond}
We now show that the DT partition function satisfies the condensation recurrence; this is now a corollary of the well-known ``graphical condensation'' theorem of Kuo:

\begin{thm}\cite[Theorem 5.1]{kuo}
\label{thm:kuo}
Let $G = (V_1, V_2, E)$ be a planar bipartite graph with a given planar embedding in which $|V_1| = |V_2|$. Let vertices $a, b, c,$ and $d$ appear in a cyclic order on a face of $G$. If $a, c \in V_1$ and $b, d \in V_2$, then
\small
\begin{equation*}
Z^{D}(G)Z^{D}(G -\{a, b, c, d\}) =Z^{D}(G  - \{a, b\})Z^{D}(G - \{c, d\})   +  Z^{D}(G  - \{a, d\})Z^{D}(G - \{b, c\}).
\end{equation*}
\normalsize
\end{thm}

Take $G$ to be $H(N; \mu_1^{rc},\mu_2^{rc},\mu_3)$ for $N$ sufficiently large\footnotemark.
Let $a$ and $b$ be the vertices in sector 1 labelled by $\max S_1^-$ and $\min S_1^+$, respectively. 
Similarly, we let $c$ and $d$ be the vertices in sector 2 labelled by $\max S_2^-$ and $\min S_2^+$. 
Then $G-\{a, b, c, d\}$ is $H(N; \mu_1,\mu_2,\mu_3)$.

The resulting six dimer-model partition functions are all instances of the topological vertex, up to order $N$.
The normalization constants $q^A, q^B, q^C$ arise because the ``folklore'' technique which associates a plane partition to a dimer configuration preserves the weight up to a factor of $q^{w(\pi_{\min})}$, where $\pi_{\min}$ is the minimal dimer configuration.
 The weight of this configuration is computed, for instance in~\cite{kuo}; for us the computation is substantially messier, as $w_{\min}$ depends on $\mu_1, \mu_2, \mu_3$ in a delicate way; we omit the details.

\footnotetext{A sufficient lower bound for $N$ depends on $\mu_1$, $\mu_2$, and $\mu_3$.}

\section{PT}
\label{sec:PT}

\subsection{Labelled $AB$ configurations}

\label{sec:PTboxconfigs}

In this section we introduce one of the main objects of our study: labelled $AB$ configurations. 

\begin{defn}
\label{conditions:ab box stacking}
 If $A \subseteq  \I^- \cup \III$ and $B \subseteq \II \cup \III$ are finite sets of boxes, then $(A, B)$ is an {\em $AB$ configuration} if the following condition is satisfied:
 \begin{itemize}
 \item[] If $w=(w_1, w_2, w_3)$ is a cell in $\I^- \cup \III$ (resp.~$w\in\II\cup\III$) and any cell in $\{ (w_1-1, w_2, w_3),
	(w_1, w_2-1, w_3), (w_1, w_2, w_3-1)\}$ supports a box in $A$ (resp.~$B$), then $w$ must support a box in $A$ (resp.~$B$).
	\end{itemize}
	
\end{defn}
We remark that this is the familiar condition for plane partitions, except that gravity is pulling the boxes in the direction $(1,1,1)$, away from the origin.

Next, we give an algorithm that labels $AB$ configurations. Note that the algorithm assigns labels to cells, not boxes.

\begin{algorithm}\label{algorithm:AB labelling algorithm}
	\begin{enumerate}
	\setlength\itemsep{-.5em}
		\item If a connected component of $(\I^-\cap A)\cup(\II\setminus B)\cup(\III\cap(A\triangle B))$ contains a box in $\Cyl_i^-\cup\II_{\bar{i}}$ and a box in $\Cyl_j\cup\II_{\bar{j}}$, where $i\neq j$, terminate with failure.
		\item For each connected component $C$ of $(\I^-\cap A)\cup(\II\setminus B)\cup(\III\cap(A\triangle B))$ that contains a box in $\Cyl_i^-\cup\II_{\bar{i}}$, label each element of $C$ by $i$.
		\item For each remaining connected component $C$ of $(\I^-\cap A)\cup(\II\setminus B)\cup(\III\cap(A\triangle B))$, label each element of $C$ by the same freely chosen element of $\mathbb{P}^1$.
	\end{enumerate}
\end{algorithm}

Because the algorithm may fail in Step 1, there are $AB$ configurations that cannot be labelled.
A {\em labelled $AB$ configuration} is an $AB$ configuration for which the labelling algorithm succeeds. Let $\sAB(\mu_1, \mu_2, \mu_3)$ denote the set of all labelled $AB$ configurations.

\begin{figure}[htb]
\centering
\includegraphics[width=0.75in]{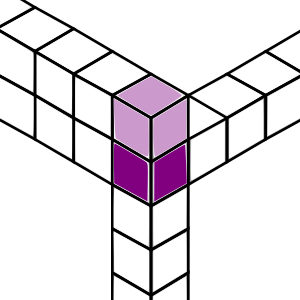}
	\hspace{0.25in}
\includegraphics[width=0.75in]{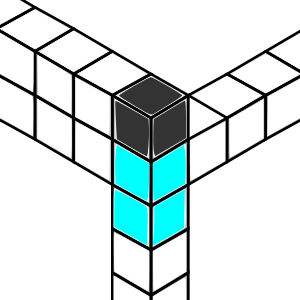}
	\hspace{0.25in}
\includegraphics[width=0.75in]{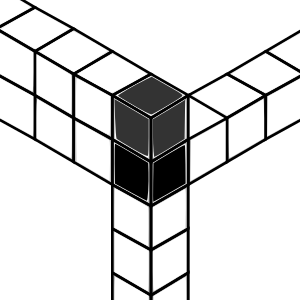}
	\hspace{0.25in}
\includegraphics[width=0.75in]{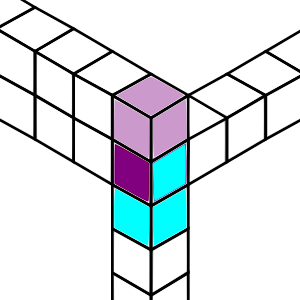}
\caption{The $AB$ configurations from Example~\ref{ex:ABlabellings}.}
\label{fig:ABlabellings}
\end{figure}

\begin{example}
\label{ex:ABlabellings}
Let $\mu_1 = (1), \mu_2 = (2),$ and $\mu_3 = (1)$. Then $\III = \{ (0, 0, 0)\}$ and $\II = \II_{\bar{1}} = \{ (0, 0, 1)\}$. In Figure~\ref{fig:ABlabellings} we illustrate four $AB$ configurations\footnotemark, three of which are labelled $AB$ configurations. 
\footnotetext{The first three of these configurations appear in \cite[Section 5.4]{PT2} as the length 1 configuration (i), the length 3 configuration (iv), and the length 2 configuration (iii).}

\begin{enumerate}
\setlength\itemsep{-.25em}
\item 
$A$ consists of a single box at $(0, 0, 0)$ and $B = \varnothing$. Step 2 of Algorithm~\ref{algorithm:AB labelling algorithm} gives the cells $(0, 0, 0)$ and $(0, 0, 1)$ the label 1, which is indicated by the color purple. The cell $(0, 0, 0)$ is opaque because it supports a box; the cell $(0, 0, 1)$ does not. 
\item 
$A = \{ (0, 0, 0), (0, 0, -1)\}$ and $B = \{ (0, 0, 1)\}$. Step 2 labels the cells in $A$ by 3, which we illustrate by coloring the two boxes cyan. The box at $(0, 0, 1)$ is colored gray because it does not get a label.
\item 
$A = \varnothing$ and $B =\{ (0, 0, 0), (0, 0, 1)\}$. Again, the box at $(0, 0, 1)$ does not get a label. The box at $(0, 0, 0)$ has a free choice of label in $\mathbb{P}^1$.
\item 
$B = \varnothing$ and $A = \{ (0, 0, 0), (0, 0, -1)\}$. The algorithm terminates with failure in Step 1 because $(0, 0, -1) \in\Cyl_{3}^{-}$ and $(0, 0, 1) \in \II_{\bar{1}}$. In the figure, $(0, 0, 0)$ is colored both cyan, required by the box at $(0, 0, -1)$, and purple, required by the cell at $(0, 0, 1)$. 
\end{enumerate}

\end{example}

 \begin{wrapfigure}{r}{.5\textwidth}
\centering
\includegraphics[width=1.5in]{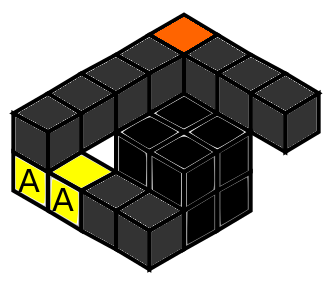}
\includegraphics[width=1.5in]{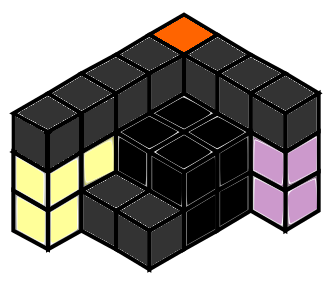}
\caption{The $AB$ configuration from Example 2}
\label{fig:tripartiteex3}
\end{wrapfigure}

\noindent{\bf Example 2.}
 Figure~\ref{fig:tripartiteex3} shows a labelled $AB$ configuration with $\mu_1 = (3, 3, 1), \mu_2 = (3, 2, 2, 1),$ and $\mu_3 = (5, 3, 3, 1)$. The left image shows the configuration. The boxes belonging to $A$ are marked; all other boxes are in $B$. The right image includes surrounding cells in $\II$. In both images, yellow cells are labelled 2 and purple cells are labelled 1. Opaque cells support a box in the configuration and transparent cells do not. The two connected components labelled by a freely chosen element of $\mathbb{P}^1$ are colored black and orange, respectively. 




Define $$W(\mu_1, \mu_2, \mu_3)=\sum_{(A, B)\in\sAB(\mu_1, \mu_2, \mu_3)}q^{\lvert A\rvert+\lvert B\rvert}.$$
We prove, in a paper in preparation, that labelled $AB$ configurations are a discrete version\footnotemark~of {\em labelled box configurations} as defined in \cite[Section 2.5]{PT2}, and therefore $W(\mu_1, \mu_2, \mu_3)$ is the {\em topological vertex in PT theory.}



\footnotetext{More precisely, there is a surjection $\psi$ from labelled $AB$ configurations to labelled box configurations, and if $\pi$ is a labelled box configuration, $|\psi^{-1}(\pi)| = \chi_{\text{top}}(\pi)$, where $\chi_{\text{top}}(\pi)$ is the topological Euler characteristic of the moduli space of labellings of $\pi$, in the terminology of \cite{PT2}. }

\subsection{PT theory and the labelled double-dimer model}



Next we explain the relationship between labelled $AB$ configurations and {\em double-dimer configurations}. On an infinite graph, a double-dimer configuration is the union of two dimer configurations. 

\begin{figure}[htb]
\centering
         \includegraphics[width=1.5in]{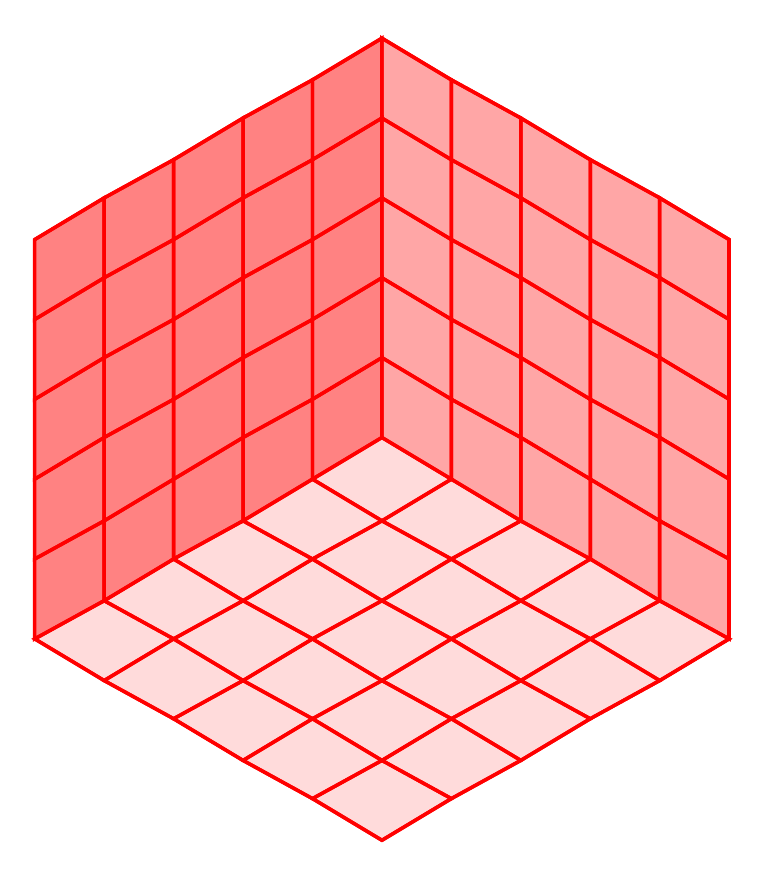}
       \includegraphics[width=1.5in]{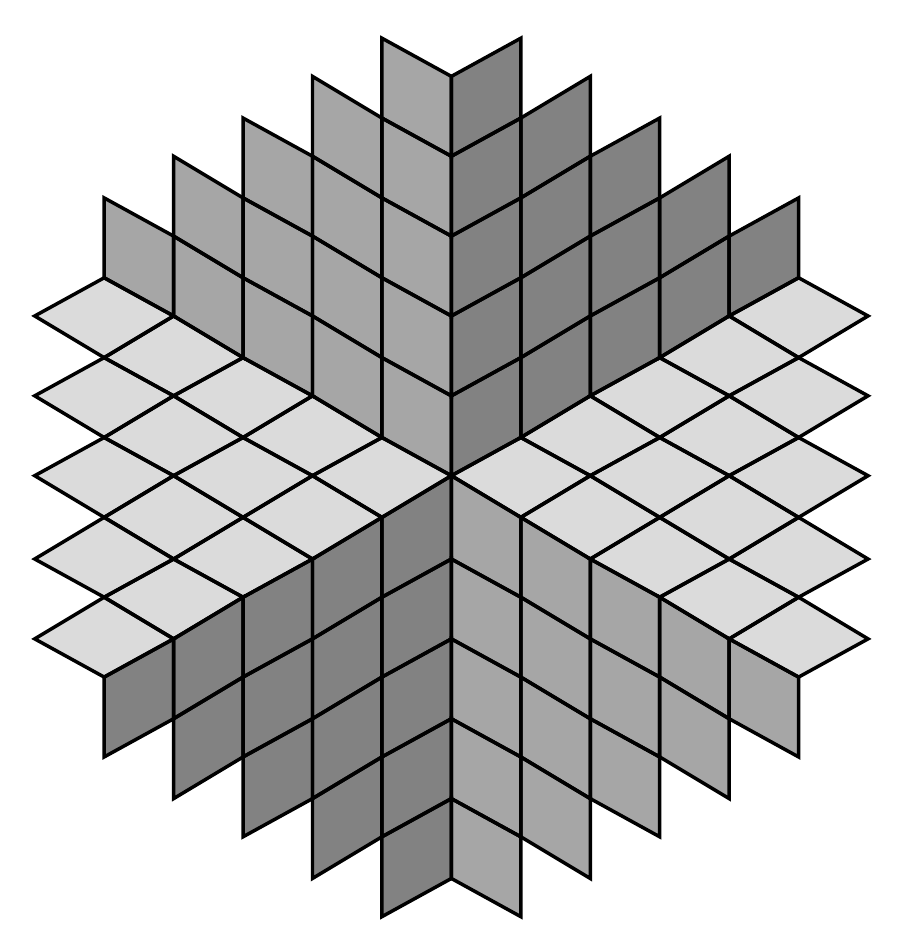}
   \includegraphics[width=1.25in]{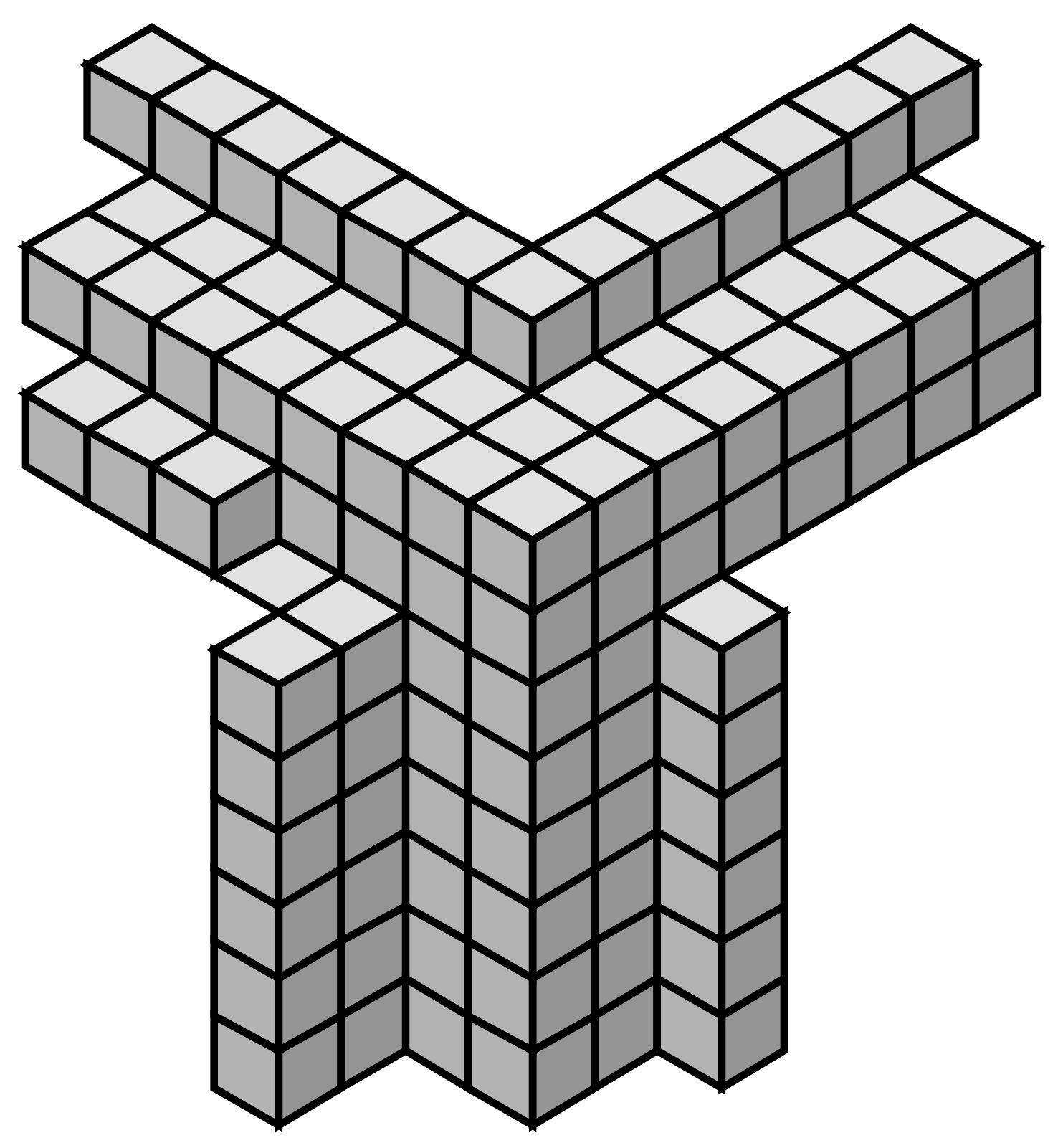}
\includegraphics[width=1.5in]{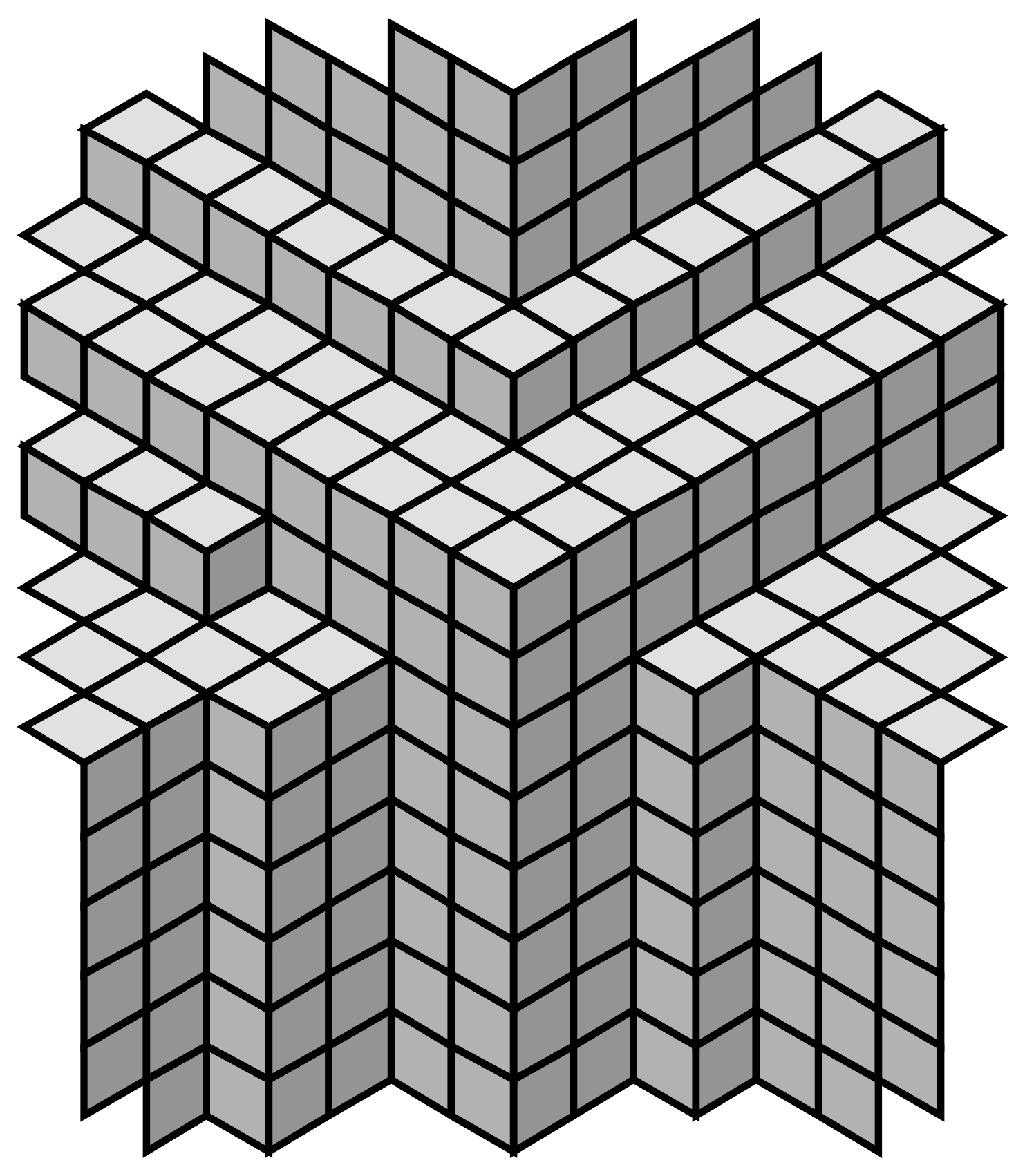}
\caption{Converting an $AB$ configuration to a tiling of the plane. Left two pictures: extra tiles that may be used to extend tilings for $B$ and $A$, respectively. Right two pictures: an example of the surface $(\I^-\cup\III)\setminus A$, and its extended tiling.}
\label{fig:empty tilings}
\end{figure}

Let $(A, B)$ be an $AB$ configuration. We consider $A$ and $B$ separately. 
For $A$, we view the surface $(\I^-\cup\III)\setminus A$ as a lozenge tiling.
In other words, we take the set of boxes $A \subseteq \I^-\cup\III$ and draw the tiles corresponding to cells that are {\underline{not}} in $A$. Similarly, for $B$, we view the surface $(\II\cup\III)\setminus B$ as a lozenge tiling. 
We then extend each of these tilings to tilings of 
the entire plane. That is, in Figure~\ref{fig:empty tilings}, we overlay the third image on the second image to obtain the final image. %
Then, these lozenge tilings 
are equivalent to dimer configurations 
of the infinite honeycomb graph $H$.



Let $M_{A}$ (resp.~$M_{B}$)  denote the dimer configuration of $H$ corresponding to the infinite tiling obtained from $A$ (resp. $B$). Superimposing $M_{A}$ and $M_{B}$ so that Region $\III$ is in the same place in the two pictures produces a double-dimer configuration $D_{(A, B)}$ on $H$.

For example, the third image of Figure~\ref{fig:empty tilings} shows the tiling corresponding to $A$, where $(A, B)$ is the $AB$ configuration from Figure~\ref{fig:tripartiteex3}. The set $A$ consists of two boxes, $(3, -1, 0)$ and $(3, 0, 0)$, so we have drawn the tiles corresponding to $(\I^{-} \cup \III) \setminus \{(3, -1, 0), (3, 0, 0)\}$.



\begin{figure}[hbt]
\centering
\includegraphics[width=1.5in]{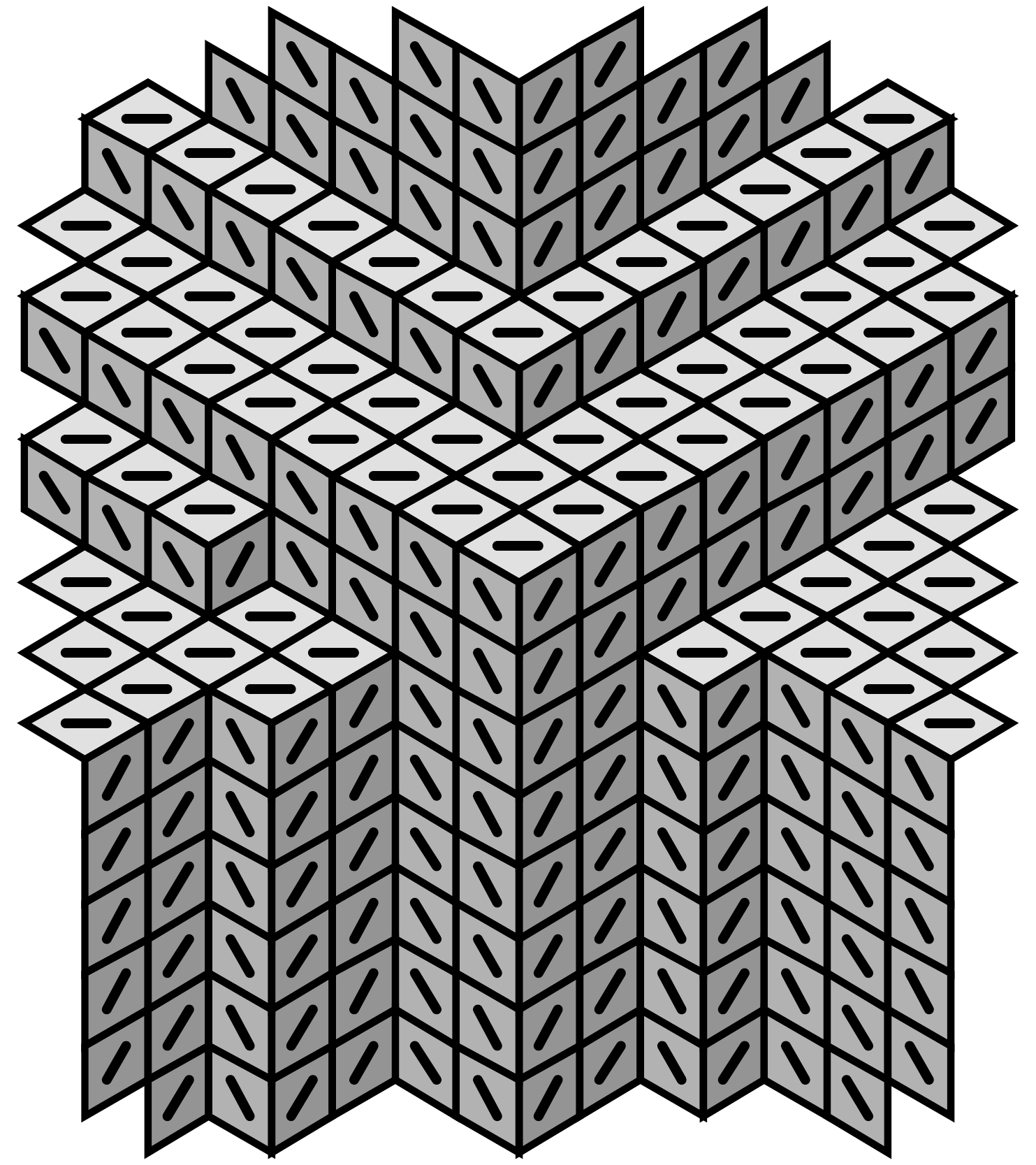}
\includegraphics[width=1.5in]{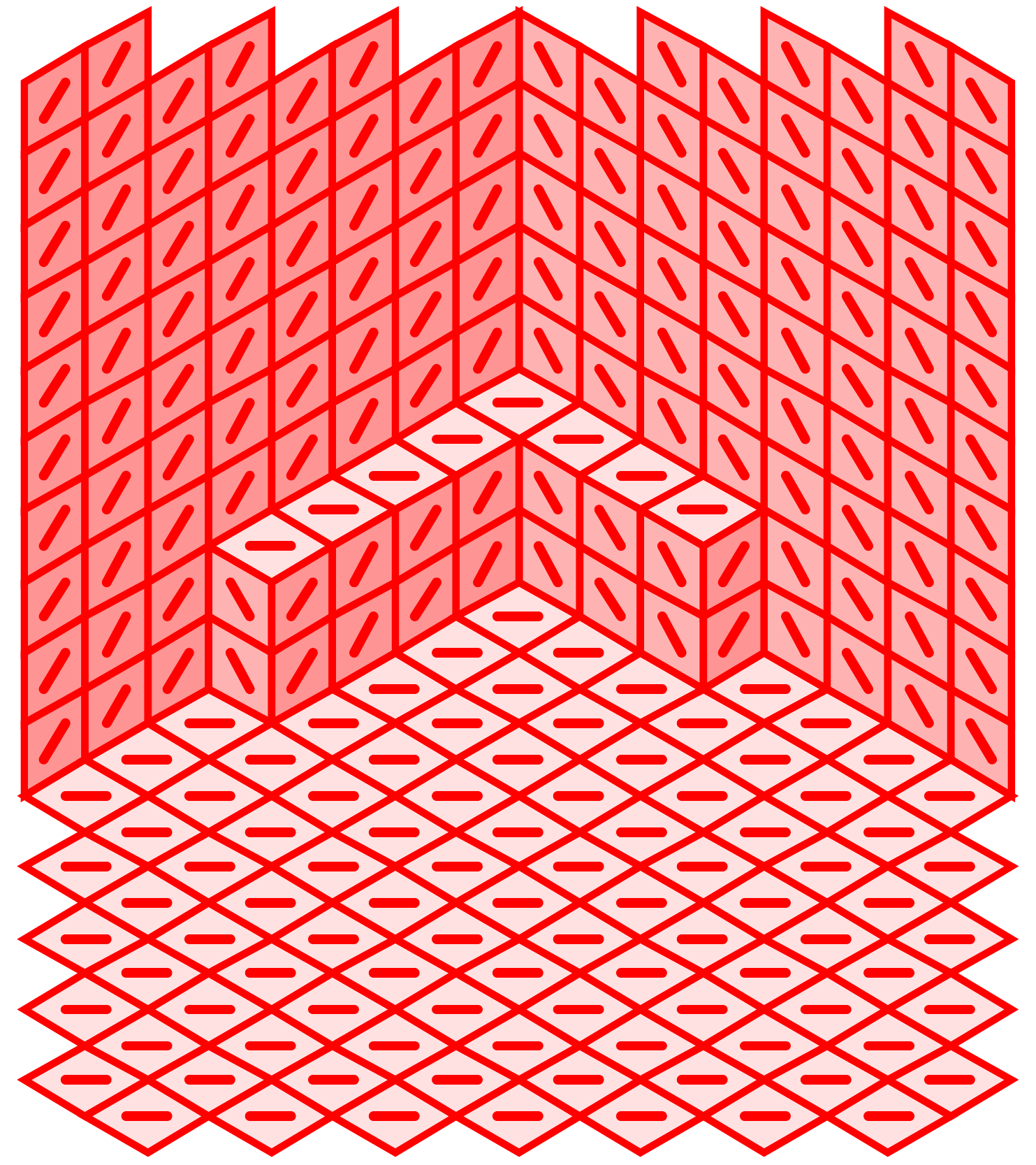}
\includegraphics[width=1.5in]{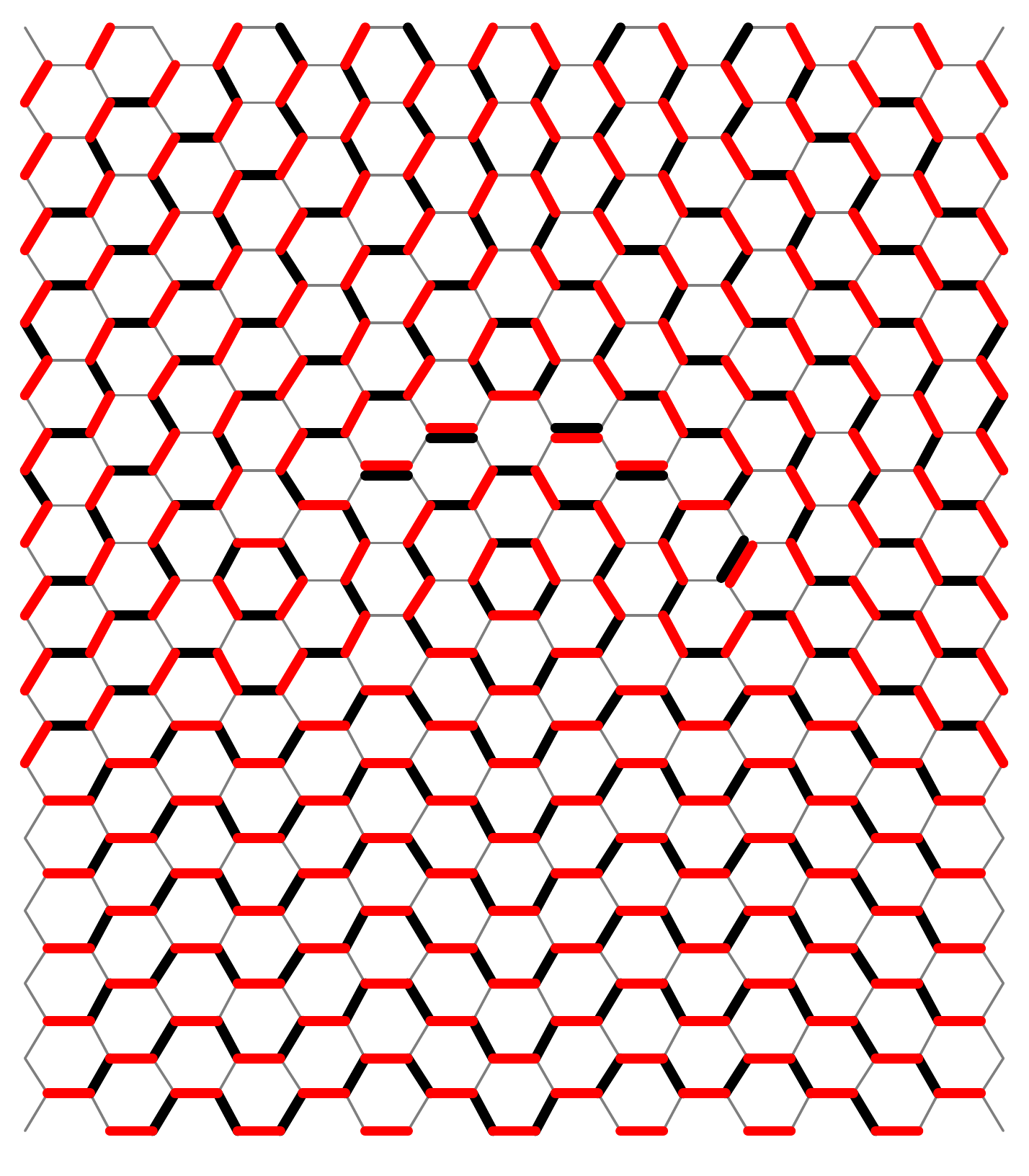}
\includegraphics[width=1.5in]{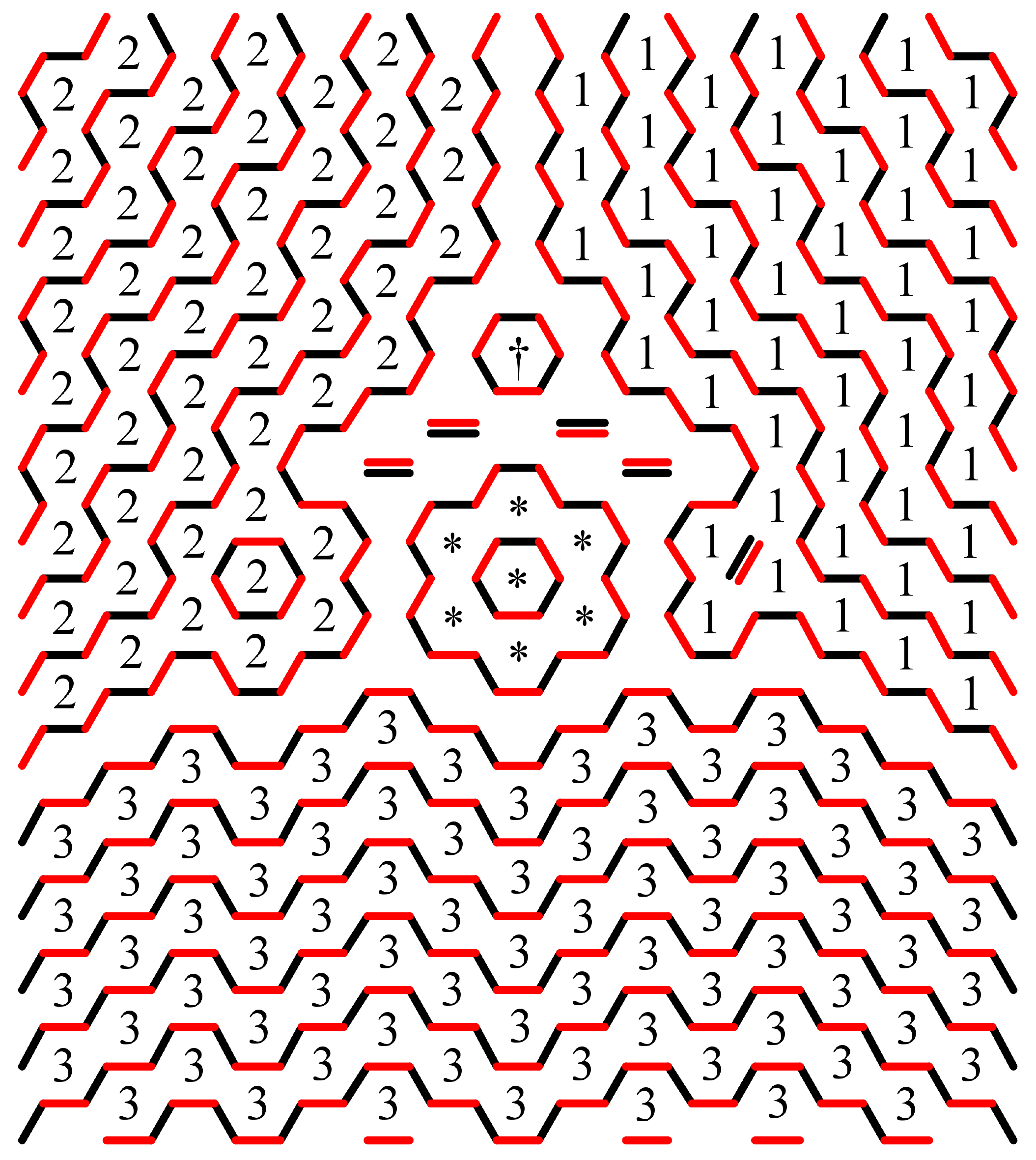}
\caption{First: The dimer configuration $M_{A}$. Second: The dimer configuration $M_{B}$. Third: The superposition of $M_A$ and $M_B$, a double-dimer configuration on $H$. Fourth: The labelled double-dimer configuration.}
\label{fig:superposition}
\end{figure}

The corresponding dimer configurations $M_{A}$ and $M_{B}$ are shown in Figure~\ref{fig:superposition}. Their superposition, shown immediately to their right, is a double-dimer configuration on $H$. 

Just as we label certain $AB$ configurations, we label certain double-dimer configurations. 
Before presenting our double-dimer labelling algorithm, we make a few remarks. 
It can be shown that each path in $D_{(A, B)}$ crosses each coordinate axis finitely many times. Consequently, 
there is a well-defined notion of the sectors\footnotemark~that contain the ends of such a path.
Also, let $h_{(A, B)}$ be the relative height function that 
assigns to each face $f$ of $H$ the height difference at $f$ of the surface corresponding to $M_B$ above that corresponding to $M_A$. The loops and paths in $D_{(A, B)}$ are the contour lines for $h_{(A, B)}$. Every path in $D_{(A, B)}$ divides the plane into two disjoint regions, and we call such a region the \emph{higher side} of the path if $h_{(A, B)}$ increases by $1$ when entering that region by crossing the path.

\footnotetext{When we refer to ``sectors'' in this section, we mean the sectors defined in the right hand side of Figure~\ref{fig:sectors}.}


\begin{algorithm}\label{algorithm:DD labelling algorithm}

	\begin{enumerate}
	\setlength\itemsep{-.5em}
		\item If $D_{(A, B)}$ contains a path whose ends are contained in different sectors, terminate with failure.
		\item For each path in $D_{(A, B)}$ such that sector $i$ contains the ends of that path, label each face of $H$ contained in the higher side of that path by $i$.
		\item For each loop in $D_{(A, B)}$ that is not contained in the interior of another loop or the higher side of a path, label each enclosed face of $H$ by the same freely chosen element of $\mathbb{P}^1$.
	\end{enumerate}
\end{algorithm}

For example, if we label the double-dimer configuration from Figure~\ref{fig:superposition}, we obtain the labelled double-dimer configuration shown in Figure~\ref{fig:superposition}. 
Observe that the paths in the double-dimer configuration in Figure~\ref{fig:superposition} are ``rainbow-like.'' In other words, the paths are nested and start and end in the same sector. 

\begin{thm}
\label{thm:tripartiteABconnection}

Let $(A, B)$ be an $AB$ configuration. 
 Then $(A, B)$ is a labelled $AB$ configuration if and only if the double-dimer configuration $D_{(A, B)}$ has the property that each path starts and ends in the same sector. 

\end{thm}

We prove Theorem~\ref{thm:tripartiteABconnection} by proving that Algorithm~\ref{algorithm:DD labelling algorithm} succeeds if and only if Algorithm~\ref{algorithm:AB labelling algorithm} succeeds.


In order to apply the double-dimer analogue of Kuo's graphical condensation,
we must truncate our double-dimer configuration on the infinite honeycomb graph to obtain 
 a {\em double-dimer configuration with nodes} on the $N \times N \times N$ honeycomb graph $H(N)$.

\begin{defn}
Let $G = (V_1, V_2, E)$ be a finite edge-weighted bipartite planar graph embedded in the plane with $|V_1| = |V_2|$. Let ${\bf N}$ denote a set of special vertices called {\em nodes} on the outer face of $G$.
A {\em double-dimer configuration} on $(G, {\bf N})$ is a multiset of the edges of $G$ with the property that each internal vertex is the endpoint of exactly two edges, and each vertex in ${\bf N}$ is the endpoint of exactly one edge.
\end{defn}


Each double-dimer configuration is associated with a planar pairing of the nodes. On a finite graph, the notion that the paths are ``rainbow-like'' means that the pairing is {\em tripartite}.

\begin{defn}
A planar pairing $\sigma$ is {\em tripartite} if the nodes can be divided into three circularly contiguous sets $R, G$, and $B$ so that no node is paired with a node in the same set. We often color the
nodes in the sets red, green, and blue, in
which case $\sigma$ is the unique planar pairing
in which like colors are not paired.
\end{defn}

\begin{wrapfigure}{l}{.25\textwidth}
\centering
\includegraphics[width=1.5in]{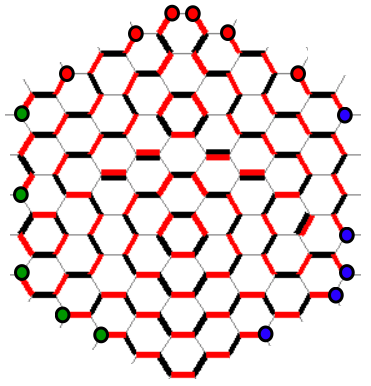}
\caption{}
\label{fig:tripartiteex1}
\end{wrapfigure}

 The process of truncating the infinite double-dimer configuration is straightforward and details are omitted here. 
 Continuing our example, truncating the double-dimer configuration from Figure~\ref{fig:superposition} to a double-dimer configuration on $H(5)$ produces the tripartite double-dimer configuration shown in Figure~\ref{fig:tripartiteex1}.

The set of nodes {\bf N} and the coloring of these nodes is determined by the partitions $\mu_1, \mu_2$, and $\mu_3$, and this can be made explicit by using the Maya diagram associated to each partition. 


We refer to the labelling and sectors of the graph $H(N)$ shown in Figure~\ref{fig:sectors}. 
To determine the nodes in Sector $i$, we draw the Maya diagram associated to $\mu_i$. 
\begin{itemize}
\setlength\itemsep{-.25em}
\item In Sector 1, the blue nodes are the holes with positive coordinates and the red nodes are the beads with negative coordinates.
\item In Sector 2, the red nodes are the holes with positive coordinates and the green nodes are the beads with negative coordinates. 
\item In Sector 3, the green nodes are the holes with positive coordinates and the blue nodes are the beads with negative coordinates. 
\end{itemize}

\subsection{The condensation identity for PT invariants}

Let $Z^{DD}_{\sigma}(G, {\bf N})$ denote the weighted sum of all double-dimer configurations with a particular pairing $\sigma$. In \cite{jenne}, the first author showed that when $\sigma$ is tripartite and certain other technical conditions hold we have the following:
\begin{eqnarray*}
& & 
Z^{DD}_{\sigma}(G, {\bf N}) Z^{DD}_{\sigma_5}(G , {\bf N} - \{x, y, w, v\})\\ & =& 
Z^{DD}_{\sigma_1}(G , {\bf N} - \{x, y\})  
Z^{DD}_{\sigma_2}(G , {\bf N} - \{w, v\})  +
 Z^{DD}_{\sigma_3}(G, {\bf N} - \{x, v\})  Z^{DD}_{\sigma_4}(G , {\bf N} - \{w, y\}) 
 \end{eqnarray*}

 

 %
  We apply this recurrence by {\em adding}
nodes to the graph so that ${\bf N} = {\bf \widetilde{N}} - \{a, b, c, d\}$ for four nodes $a$, $b$, $c$ and $d$. We choose the four nodes as follows:
  Let $a$ and $b$ be the nodes in sector 1 labelled by $\max S_1^-$ and $\min S_1^+$, respectively. 
Similarly, we let $c$ and $d$ be the nodes in sector 2 labelled by $\max S_2^-$ and $\min S_2^+$. Note that these nodes have the same coordinates as the vertices specified in Section~\ref{sec:DTcond} but the coordinate system is different (see Figure~\ref{fig:sectors}). Many details here have been omitted, due to space constraints.




 
 





\acknowledgements{We would like to thank Jim Bryan, Rick Kenyon, Rahul Pandharipande, Richard Thomas, Jim Propp, Karel Faber, Kurt Johansson, and frankly countless other geometers, probabilists and combinatorialists for helpful conversations.

\bibliographystyle{alpha}
\bibliography{ptdt} 

\end{document}